\documentclass[11pt]{article}
\usepackage{graphicx}
\usepackage{amssymb}
\usepackage{fancyhdr}
\newtheorem{theorem}{Theorem}[section]
\newtheorem{definition}[theorem]{Definition}
\newtheorem{proposition}[theorem]{Proposition}
\newtheorem{lemma}[theorem]{Lemma}

\begin{document}

\title{Differential Inequalities for Distance Comparison}

\author{\c{S}ahin Ko\c{c}ak\thanks{Department of Mathematics, Anadolu University
26470 Eski\c{s}ehir, Turkey\newline\indent\,\,\, E-mail: skocak@anadolu.edu.tr}
\quad Murat Limoncu\thanks{Department of Mathematics, Anadolu University
26470 Eski\c{s}ehir, Turkey\newline\indent\,\,\, E-mail: mlimoncu@anadolu.edu.tr}}

\maketitle

\begin{abstract}
\noindent Comparison of $1$-dimensional distance functions is a basic tool in  Alexandrov geometry
and it is used to characterize spaces with curvature bounded above or below. For the zero curvature bound
there is a differential inequality which enables one to check this comparison directly on a given smooth
$1$-dimensional distance function. In this note we give a generalization of this property to arbitrary curvature bounds.
\end{abstract}


\vskip3mm


\section{Introduction}
Alexandrov geometry is based on distance comparison.
Let $(X,d)$ be a metric space with an intrinsic metric i.e. distance between two points is realized as the length of a path between these points.
(These spaces are also called geodesic spaces. For the basic notions of metric geometry see \cite{BBI} and \cite{Papadopoulos}).
The prime example of geodesic spaces is a complete Riemannian manifold with the metric induced by the Riemannian structure.

A segment in $X$ is a shortest path between its end points, parametrized by arc length
(a segment with end points $x,y\in X$ is denoted by $[xy]$, which however might not be determined uniquely by $x$ and $y$).
Given a point $p\in X$ and a segment $\gamma:[a,b]\longrightarrow X$ with end points $x=\gamma(a)$ and $y=\gamma(b)$, we consider
the function $g(t)=d(p,\gamma(t))$. One of the tools of Alexandrov geometry is to compare this ``{\it 1-dimensional distance function}"
with an appropriate Euclidean 1-dimensional distance function.
To do so, a comparison segment $\gamma_0=[\bar{x}\bar{y}]$ of the same length as $\gamma$ is chosen in the Euclidean plane, and a reference point $\bar{p}$,
which is positioned in the same way to $\gamma_0$ as $p$ to $\gamma$, i.e. $d_0(\bar{p},\bar{x})=d(p,x)$ and $d_0(\bar{p},\bar{y})=d(p,y)$.
($d_0$ denotes the Euclidean distance) The function $g_0:[a,b]\longrightarrow\mathbb{R}$, $g_0(t)=d_0(\bar{p},\gamma_0(t))$
is called the comparison function for $g$.

The following is one of the basic characterizations of nonpositive (resp. nonnegative) curvature in Alexandrov geometry (see \cite{BBI}).
\begin{definition}\label{defofRk=0}
  A geodesic space $X$ is called nonpositively (resp. nonnegatively) curved, if every point in $X$ has a neighborhood such that, whenever
  a point $p$ and a segment $\gamma$ lie within this neighborhood, the comparison function $g_0$ for $g$ ($g(t)=d(p,\gamma(t))$)
  satisfies $g_0(t)\geq g(t)$ (resp. $g_0(t)\leq g(t)$) for all $t\in[a,b]$.
\end{definition}

One can choose instead of the Euclidean plane other two dimensional simply connected space forms for comparison.
Let $\mathbb{H}_k$ denote the hyperbolic plane with curvature $k<0$ (i.e. metric tensor of the standard hyperbolic plane multiplied by $-\frac{1}{k}$) and
$\mathbb{S}_k$ denote a Euclidean sphere with curvature $k>0$ (i.e. radius $\frac{1}{\sqrt{k}}$) with its intrinsic metric.
We will denote these model spaces by
$\mathbb{M}_k$ ($\mathbb{M}_k$=$\mathbb{R}^2$ for $k=0$, $\mathbb{M}_k$=$\mathbb{H}_k$ for $k<0$ and $\mathbb{M}_k$=$\mathbb{S}_k$ for $k>0$)
and the intrinsic metric on $\mathbb{M}_k$ by $d_k$.

Given a segment $\gamma$ in $X$ and a point $p\in X$, one can construct a comparison function $g_k$ for $\mathbb{M}_k$ as in the Euclidean case.
(One can always find a comparison segment and a reference point in $\mathbb{M}_k$ for $k\leq0$;
in $\mathbb{S}_k$, if the segment $\gamma$ is small enough and $p$ close enough to $\gamma$.) One can then generalize the above definition
verbatim to the model space $\mathbb{M}_k$:
\begin{definition}\label{defofRk=arbitrary}
  A geodesic space $X$ is called of ${\rm curvature}\leq k$ (resp. of ${\rm curvature}\geq k$), if every point in $X$ has a neighborhood such that, whenever
  a point $p$ and a segment $\gamma$ lie within this neighborhood, the comparison function $g_k$ in $\mathbb{M}_k$ satisfies
  $g_k(t)\geq g(t)$ (resp. $g_k(t)\leq g(t)$) for all $t\in[a,b]$.
\end{definition}

\begin{figure}[h!]
 \begin{center}
  \includegraphics[width=0.75\textwidth]{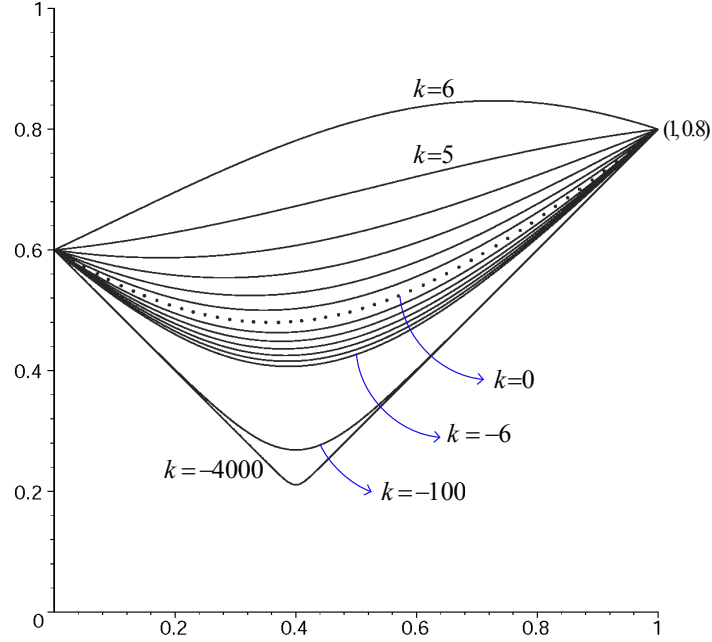}
 \end{center}
 \caption{Comparison functions $g_k$ in model spaces $\mathbb{M}_k$ for a segment of unit length and with $g_k(0)=3/5$, $g_k(1)=4/5$.
 The sample values of $k$ are $6,5,\dots,1,0,-1,\dots,-5,-6,-100,-4000$.}
 \label{figure}
\end{figure}

We note that the comparison function $g_k$ associated to $g$ in $\mathbb{M}_k$ is completely determined by three numbers
$b-a$, $d(p,\gamma(a))$ and $d(p,\gamma(b))$ (under the clause referred to above for $\mathbb{S}_k$).
An example is given in Figure \ref{figure} for $b-a=1$, $d(p,\gamma(a))=3/5$, $d(p,\gamma(b))=4/5$
and for some $k$. (The collection of the comparison functions constitute a kind of curvature scale.)

We can express the comparison functions on the model spaces as follows:
\begin{eqnarray}
  g_0(t)=\sqrt{(u-t)^2+v^2}\qquad{\rm for}\,\, k=0\,\,\,{\rm (on}\,\,\mathbb{R}^2{\rm )},
\end{eqnarray}
\begin{eqnarray}
  \!\!\!\!\!\!\!\!\!\!g_k(t)\!\!\!\!&=&\!\!\!\!\frac{1}{\sqrt{-k}}{\rm argcosh}\Big(\frac{1}{2v}\Big[(u^2+v^2)e^{-\sqrt{-k}t}+e^{\sqrt{-k}t}\Big]\Big)\,\,\,{\rm for}\,\,k<0\,\,\,{\rm (on}\,\,\mathbb{H}_k{\rm )},\\
  \!\!\!\!\!\!\!\!\!\!g_k(t)\!\!\!\!&=&\!\!\!\!\frac{1}{\sqrt{k}}\arccos\Big(\sqrt{k}\Big[u\cos(\sqrt{k}t)+v\sin(\sqrt{k}t)\Big]\Big)\,\,\,{\rm for}\,\, k>0\,\,\,{\rm (on}\,\,\mathbb{S}_k{\rm )},
\end{eqnarray}
whereby $u$ and $v$ are parameters to adjust the values $g_k(a)=d(p,\gamma(a))$ and $g_k(b)=d(p,\gamma(b))$
for a given segment $\gamma$ in $X$ and the point $p\in X$.
(We restrict the parameter values by $v>0$ for $g_k$ with $k\leq0$ and by $u^2+v^2<\frac{1}{k}$ for $g_k$ with $k>0$.)

To obtain these expressions we have made the following comparison choices:

1. On $\mathbb{R}^2$:
Comparison segment $\gamma_0:[a,b]\longrightarrow \mathbb{R}^2$, $\gamma_0(t)=(t,0)$.
Comparison point $(u,v)$ with $v>0$.
\vskip0.3cm
2. On $\mathbb{H}_k$:
Comparison segment $\gamma_k:[a,b]\longrightarrow \mathbb{H}_k$, $\gamma_k(t)=(0,e^{\sqrt{-k}t})$.
Comparison point $(u,v)$ with $v>0$.
\vskip0.3cm
3. On $\mathbb{S}_k$:
Comparison segment $\gamma_k:[a,b]\longrightarrow \mathbb{S}_k$ (with $b-a\leq\frac{\pi}{\sqrt{k}}$),
$\gamma_k(t)=(\frac{1}{\sqrt{k}}\cos\sqrt{k}t,\frac{1}{\sqrt{k}}\sin\sqrt{k}t,0)$.
Comparison point $(u,v,\sqrt{\frac{1}{k}-u^2-v^2})$.
(We assume $a\geq0$; often it will be $a=0$.)
Under these choices, $g_k(t)$ is the distance of the comparison point to $\gamma_k(t)$ (with respect to $d_k$).

We remark that the functions $g_k$ (for any $k$) are, by the triangle inequality, nonexpanding in the sense that $|g_k(t_1)-g_k(t_2)|\leq|t_1-t_2|$
holds for any $t_1,t_2\in[a,b]$. By the same reason and since the comparison segment $\gamma_k$ is arc-length parametrized, we also have
\begin{eqnarray}
  |t_1-t_2|\leq g_k(t_1)+g_k(t_2).
\end{eqnarray}

This last property is not a consequence of nonexpandingness as, for example, the simple function
\begin{eqnarray}
  g:[0,1]\longrightarrow\mathbb{R},\quad g(t)=\frac{1}{4}(t+1)
\end{eqnarray}
shows. From the point of view of the theorems we give in the next section, we want to consider the following type of functions:
\begin{definition}
  A non-negative function $g:[a,b]\longrightarrow\mathbb{R}$ is called distance-like, if it is nonexpanding and the following property is satisfied:
\vskip0.2cm
  For any $t_1,t_2\in[a,b]$, it holds $|t_1-t_2|\leq g(t_1)+g(t_2)$.
\end{definition}
In other words, the function $g$ is distance-like, if for any $a\leq t_1\leq t_2\leq b$, the three numbers $g(t_1)$, $g(t_2)$ and $t_2-t_1$
can be realized as the edge-lengths of an Euclidean triangle. (One can then find a triangle in $\mathbb{H}_k$ with the same edge-lengths for any $k<0$;
and a triangle in $\mathbb{S}_k$ if these numbers are small enough, depending on $k>0$.)

One can simplify the above condition by the following simple lemma:
\begin{lemma}
  Let $g:[a,b]\longrightarrow\mathbb{R}$ be a non-negative, nonexpanding function satisfying $b-a\leq g(a)+g(b)$.
  Then $g$ is a distance-like function.
\end{lemma}
{\it Proof.} Let $a\leq t_1\leq t_2\leq b$. By nonexpandingness we can write
\begin{eqnarray}
  |g(t_1)-g(a)|\leq t_1-a\quad{\rm and}\quad|g(b)-g(t_2)|\leq b-t_2.
\end{eqnarray}
Adding up the inequalities $g(t_1)\geq g(a)-t_1+a$ and $g(t_2)\geq g(b)+t_2-b$, we get $g(t_1)+g(t_2)\geq t_2-t_1$.

By adding an appropriate positive constant, one can obviously make a nonexpanding function distance-like.

\section{Differential Inequalities for Function Comparison}
In \cite{BBI}, the following property is given (in a slightly different formulation) as an interesting exercise:
\begin{proposition}\label{prpk=0}
  1. The function $g_0(t)=\sqrt{(u-t)^2+v^2}$ satisfies the differential equation
  \begin{eqnarray}
    g''_0=\frac{1-(g'_0)^2}{g_0}.
  \end{eqnarray}
  2. Let $g:[a,b]\longrightarrow\mathbb{R}$ be a smooth (or of class $\mathcal{C}^2$),
  positive-valued, distance-like function.

  i) The differential inequality
  \begin{eqnarray}
    g''\geq\frac{1-(g')^2}{g}
  \end{eqnarray}
  is satisfied on $[a,b]$ if and only if the following holds: Given any $a\leq t_1<t_2\leq b$, then $g(t)\leq g_0(t)$
  on $[t_1,t_2]$, whereby the parameter values of $g_0$ are so adjusted that $g_0(t_1)=g(t_1)$ and $g_0(t_2)=g(t_2)$.

  ii) The differential inequality
  \begin{eqnarray}
    g''\leq\frac{1-(g')^2}{g}
  \end{eqnarray}
  is satisfied on $[a,b]$ if and only if the following holds: Given any $a\leq t_1<t_2\leq b$, then $g(t)\geq g_0(t)$
  on $[t_1,t_2]$, whereby the parameter values of $g_0$ are so adjusted that $g_0(t_1)=g(t_1)$ and $g_0(t_2)=g(t_2)$.
\end{proposition}

We want to generalize this property to the other function types $g_k$.
For definiteness we want to introduce the following notation:
Given a smooth, positive valued distance-like function $g:[a,b]\longrightarrow\mathbb{R}$ ($0\leq a<b$) and $a\leq t_1<t_2\leq b$,
we denote the function $g_k$ with parameter values adjusted to $g$ on $t_1$ and $t_2$, i.e. $g_k(t_1)=g(t_1)$ and $g_k(t_2)=g(t_2)$, by $g^{t_1,t_2}_k$.
For $k\leq0$ this is always possible; for $k>0$, $b$ and $\displaystyle\max_{a\leq t\leq b}g(t)$ must be small enough.

We first consider the hyperbolic case:

\begin{proposition}\label{prpk<0}
  Let $k<0$. Then,

  \noindent
  1. The function $g_k$ satisfies the differential equation
  \begin{eqnarray}
    g''_k=\sqrt{-k}\Big(1-(g'_k)^2\Big)\coth(\sqrt{-k}g_k).
  \end{eqnarray}
  2. Let $g:[a,b]\longrightarrow\mathbb{R}$ be a smooth (or of class $\mathcal{C}^2$),
  positive-valued, distance-like function.

  i) The differential inequality
  \begin{eqnarray}\label{inprp22kmns}
    g''\geq\sqrt{-k}\Big(1-(g')^2\Big)\coth(\sqrt{-k}g)
  \end{eqnarray}
  is satisfied on $[a,b]$ if and only if $g\leq g^{t_1,t_2}_k$ holds on $[t_1,t_2]$ for any $a\leq t_1<t_2\leq b$.

  ii) The differential inequality
  \begin{eqnarray}
    g''\leq\sqrt{-k}\Big(1-(g')^2\Big)\coth(\sqrt{-k}g)
  \end{eqnarray}
  is satisfied on $[a,b]$ if and only if $g\geq g^{t_1,t_2}_k$ holds on $[t_1,t_2]$ for any $a\leq t_1<t_2\leq b$.
\end{proposition}
{\it Proof.} The first part is straightforward.
For the proof of the case (i) in part 2 of the proposition first assume that $g(t)\leq g^{t_1,t_2}_k(t)$ holds on $[t_1,t_2]$.
Then one can write
\begin{eqnarray}\label{frml1and2kmns}
  g'(t_1)\leq (g^{t_1,t_2}_k)'(t_1)\quad{\rm and}\quad g'(t_2)\geq (g^{t_1,t_2}_k)'(t_2).
\end{eqnarray}
By definition, the function $g^{t_1,t_2}_k$ satisfies
\begin{eqnarray}\label{combcoskmns}
  \cosh(\sqrt{-k}g^{t_1,t_2}_k(t))=\frac{1}{2v}\Big((u^2+v^2)e^{-\sqrt{-k}t}+e^{\sqrt{-k}t}\Big)
\end{eqnarray}
and by derivation
\begin{eqnarray}\label{combdersinkmns}
  (g^{t_1,t_2}_k)'(t)\sinh(\sqrt{-k}g^{t_1,t_2}_k(t))=\frac{1}{2v}\Big(-(u^2+v^2)e^{-\sqrt{-k}t}+e^{\sqrt{-k}t}\Big).
\end{eqnarray}
Adding up (\ref{combcoskmns}) and (\ref{combdersinkmns}), we get
\begin{eqnarray}\label{prdcskmns}
  e^{-\sqrt{-k}t}\Big(\cosh(\sqrt{-k}g^{t_1,t_2}_k(t))+(g^{t_1,t_2}_k)'(t)\sinh(\sqrt{-k}g^{t_1,t_2}_k(t))\Big)=\frac{1}{v}.
\end{eqnarray}
Since $g'(t_1)\leq (g^{t_1,t_2}_k)'(t_1)$ and $g^{t_1,t_2}_k(t_1)=g(t_1)$, we can write
\begin{eqnarray}\label{inqatt1prdcskmns}
  e^{-\sqrt{-k}t_1}\Big(\cosh(\sqrt{-k}g(t_1))+g'(t_1)\sinh(\sqrt{-k}g(t_1))\Big)\leq\frac{1}{v}.
\end{eqnarray}
Similarly, for $t=t_2$, we obtain
\begin{eqnarray}\label{odinqatt2prdcskmns}
  e^{-\sqrt{-k}t_2}\Big(\cosh(\sqrt{-k}g(t_2))+g'(t_2)\sinh(\sqrt{-k}g(t_2))\Big)\geq\frac{1}{v}.
\end{eqnarray}
The inequalities (\ref{inqatt1prdcskmns}) and (\ref{odinqatt2prdcskmns}) show that the function
\begin{eqnarray}\label{Htkmns}
  h(t):=e^{-\sqrt{-k}t}\Big(\cosh(\sqrt{-k}g(t))+g'(t)\sinh(\sqrt{-k}g(t))\Big)
\end{eqnarray}
is non-decreasing on $[a,b]$. Thus it must be
\begin{eqnarray}
  h'(t)&=&e^{-\sqrt{-k}t}\sinh(\sqrt{-k}g(t))g''(t)\nonumber\\
  &&+e^{-\sqrt{-k}t}\sqrt{-k}\cosh(\sqrt{-k}g(t))\big[(g'(t))^2-1\big]\geq0
\end{eqnarray}
on $[a,b]$. From this, we obtain, for all $t\in[a,b]$,
\begin{eqnarray}\label{orkmns}
  g''(t)\geq\sqrt{-k}\Big(1-(g'(t))^2\Big)\coth(\sqrt{-k}g(t))
\end{eqnarray}
that is (\ref{inprp22kmns}), as required.

We now consider the other direction of the case (i) of part 2 of the proposition.
If we have the differential inequality (\ref{inprp22kmns}) for all $t\in[a,b]$,
then by multiplication of both sides of (\ref{inprp22kmns})
with the positive valued function ``$e^{-\sqrt{-k}t}\sinh(\sqrt{-k}g(t))$", we get
\begin{eqnarray}
  e^{-\sqrt{-k}t}\Big(\sinh(\sqrt{-k}g(t))g''(t)+\sqrt{-k}\cosh(\sqrt{-k}g(t))[(g'(t))^2-1]\Big)\geq0.
\end{eqnarray}
The left hand side of this inequality equals to the derivative of the function $h(t)$ defined in (\ref{Htkmns}).
Hence one has $h'(t)\geq0$ and consequently, the function $h(t)$ is non-decreasing. On the other hand,
the function $h(t)$ can be written as
\begin{eqnarray}
  h(t)=\frac{\rho'(t)}{\varphi'(t)},
\end{eqnarray}
where
\begin{eqnarray}\label{psiphikmns}
  \rho(t):=e^{\sqrt{-k}t}\cosh(\sqrt{-k}g(t))\quad{\rm and}\quad\varphi(t):=\frac{1}{2}e^{2\sqrt{-k}t}.
\end{eqnarray}
Now, by using the extended mean value theorem on the intervals $[t_1,t]$ and $[t,t_2]$ for $t_1<t<t_2$, we get
\begin{eqnarray}\label{etakmns}
  \frac{\rho(t)-\rho(t_1)}{\varphi(t)-\varphi(t_1)}=\frac{\rho'(\eta)}{\varphi'(\eta)}=h(\eta),\quad{\rm for\,\,some}\quad\eta\in(t_1,t)
\end{eqnarray}
and
\begin{eqnarray}\label{xikmns}
  \frac{\rho(t_2)-\rho(t)}{\varphi(t_2)-\varphi(t)}=\frac{\rho'(\xi)}{\varphi'(\xi)}=h(\xi),\quad{\rm for\,\,some}\quad\xi\in(t,t_2).
\end{eqnarray}
Since the function $h(t)$ is non-decreasing, we have $h(\eta)\leq h(\xi)$, i.e.
\begin{eqnarray}
  \frac{\rho(t)-\rho(t_1)}{\varphi(t)-\varphi(t_1)}\leq\frac{\rho(t_2)-\rho(t)}{\varphi(t_2)-\varphi(t)}
\end{eqnarray}
for all $t\in(t_1,t_2)$, which we can also write in the form
\begin{eqnarray}\label{gvngthrsltkmns}
  \rho(t)\leq\frac{1}{\varphi(t_2)-\varphi(t_1)}\Big([\varphi(t_2)-\varphi(t)]\rho(t_1)+[\varphi(t)-\varphi(t_1)]\rho(t_2)\Big)
\end{eqnarray}
since the function $\varphi(t)=\frac{1}{2}e^{2\sqrt{-k}t}$ is strictly increasing. Using the definitions of $\rho$
and $\varphi$ in both sides of (\ref{gvngthrsltkmns}), we get
\begin{eqnarray}
  e^{\sqrt{-k}t}\cosh(\sqrt{-k}g(t))\leq\frac{2\mathcal{A}}{e^{2\sqrt{-k}t_2}-e^{2\sqrt{-k}t_1}},
\end{eqnarray}
where
\begin{eqnarray}
  \mathcal{A}&=&\frac{1}{2}e^{\sqrt{-k}t_1}\Big(e^{2\sqrt{-k}t_2}-e^{2\sqrt{-k}t})\Big)\cosh(\sqrt{-k}g(t_1))\nonumber\\
  &&+\frac{1}{2}e^{\sqrt{-k}t_2}\Big(e^{2\sqrt{-k}t}-e^{2\sqrt{-k}t_1})\Big)\cosh(\sqrt{-k}g(t_2)).
\end{eqnarray}
Here by using
\begin{eqnarray}
  \!\!\!\!\!\!\!\!\!\!\!\!\!\!\!\!&\,&g(t_1)=g^{t_1,t_2}_k(t_1)=\frac{1}{\sqrt{-k}}{\rm argcosh}\Big(\frac{1}{2v}\Big[(u^2+v^2)e^{-\sqrt{-k}t_1}+e^{\sqrt{-k}t_1}\Big]\Big)\\
  \!\!\!\!\!\!\!\!\!\!\!\!\!\!\!\!&&{\rm (and\,\,the\,\,same\,\,for}\,\,t_2{\rm )},\nonumber
\end{eqnarray}
we get
\begin{eqnarray}
  e^{\sqrt{-k}t}\cosh(\sqrt{-k}g(t))\leq\frac{1}{2v}\Big(u^2+v^2+e^{2\sqrt{-k}t}\Big),
\end{eqnarray}
whence we obtain
\begin{eqnarray}
  g(t)\leq\frac{1}{\sqrt{-k}}{\rm argcosh}\Big(\frac{1}{2v}\Big[(u^2+v^2)e^{-\sqrt{-k}t}+e^{\sqrt{-k}t}\Big]\Big)=g^{t_1,t_2}_k(t)
\end{eqnarray}
on $[t_1,t_2]$.

The proof of the case (ii) in the second part of the proposition
is proved along the same lines, only by reversing the inequalities.

We now give the corresponding theorem for $k>0$.

\begin{proposition}\label{prpk>0}
  Let $k>0$. Then,

  \noindent
  1. The function $g_k$ satisfies the differential equation
  \begin{eqnarray}
    g''_k=\sqrt{k}\Big(1-(g'_k)^2\Big)\cot(\sqrt{k}g_k).
  \end{eqnarray}
  2. Let $g:[a,b]\longrightarrow\mathbb{R}$ be a smooth (or of class $\mathcal{C}^2$),
  positive-valued, distance-like function.

  i) The differential inequality
  \begin{eqnarray}\label{inprp22kpls}
    g''\geq\sqrt{k}\Big(1-(g')^2\Big)\cot(\sqrt{k}g)
  \end{eqnarray}
  is satisfied on $[a,b]$ if and only if $g\leq g^{t_1,t_2}_k$ holds on $[t_1,t_2]$ for any $a\leq t_1<t_2\leq b$.

  ii) The differential inequality
  \begin{eqnarray}
    g''\leq\sqrt{k}\Big(1-(g')^2\Big)\cot(\sqrt{k}g)
  \end{eqnarray}
  is satisfied on $[a,b]$ if and only if $g\geq g^{t_1,t_2}_k$ holds on $[t_1,t_2]$ for any $a\leq t_1<t_2\leq b$.
\end{proposition}

(We assume here, as noted earlier, $b$ and $\displaystyle\max_{a\leq t\leq b}g(t)$ to be small enough, depending on $k$, to make the triangles
with edge lengths $t_2-t_1$, $g(t_1)$, $g(t_2)$ embeddable into the model space $\mathbb{S}_k$. We additionally want them to be small enough to make
the values $\cos(\sqrt{k}t)$ and $\sin(\sqrt{k}g(t))$ positive, as needed in the proof below.)

{\it Proof.} The first part is straightforward.
For the proof of the case (i) in part 2 of the proposition first assume that $g(t)\leq g^{t_1,t_2}_k(t)$ holds on $[t_1,t_2]$.
Then one can write
\begin{eqnarray}\label{frml1and2}
  g'(t_1)\leq (g^{t_1,t_2}_k)'(t_1)\quad{\rm and}\quad g'(t_2)\geq (g^{t_1,t_2}_k)'(t_2).
\end{eqnarray}
By definition, the function $g^{t_1,t_2}_k$ satisfies
\begin{eqnarray}\label{combcoskpls}
  \cos(\sqrt{k}g^{t_1,t_2}_k(t))=\sqrt{k}\Big(u\cos(\sqrt{k}t)+v\sin(\sqrt{k}t)\Big)
\end{eqnarray}
and by derivation
\begin{eqnarray}\label{combdersinkpls}
  -(g^{t_1,t_2}_k)'(t)\sin(\sqrt{k}g^{t_1,t_2}_k(t))=\sqrt{k}\Big(-u\sin(\sqrt{k}t)+v\cos(\sqrt{k}t)\Big).
\end{eqnarray}
Combining (\ref{combcoskpls}) and (\ref{combdersinkpls}), we get (here and below we use $\cos(\sqrt{k}t)>0$ and $\sin(\sqrt{k}g(t))>0$)
\begin{eqnarray}\label{prdcskpls}
  \!\!\!\!(g^{t_1,t_2}_k)'(t)\cos(\sqrt{k}t)\sin(\sqrt{k}g^{t_1,t_2}_k(t))-\sin(\sqrt{k}t)\cos(\sqrt{k}g^{t_1,t_2}_k(t))=-v\sqrt{k}.
\end{eqnarray}
Since $g'(t_1)\leq (g^{t_1,t_2}_k)'(t_1)$ and $g^{t_1,t_2}_k(t_1)=g(t_1)$, we can write
\begin{eqnarray}\label{inqatt1prdcskpls}
  g'(t_1)\cos(\sqrt{k}t_1)\sin(\sqrt{k}g(t_1))-\sin(\sqrt{k}t_1)\cos(\sqrt{k}g(t_1))\leq-v\sqrt{k}.
\end{eqnarray}
Similarly, for $t=t_2$, we obtain
\begin{eqnarray}\label{odinqatt2prdcskpls}
  g'(t_2)\cos(\sqrt{k}t_2)\sin(\sqrt{k}g(t_2))-\sin(\sqrt{k}t_2)\cos(\sqrt{k}g(t_2))\geq-v\sqrt{k}.
\end{eqnarray}
The inequalities (\ref{inqatt1prdcskpls}) and (\ref{odinqatt2prdcskpls}) show that the function
\begin{eqnarray}\label{Htkpls}
  H(t):=g'(t)\cos(\sqrt{k}t)\sin(\sqrt{k}g(t))-\sin(\sqrt{k}t)\cos(\sqrt{k}g(t))
\end{eqnarray}
is non-decreasing on $[a,b]$. Thus it must be
\begin{eqnarray}
  \!\!\!\!H'(t)=\cos(\sqrt{k}t)\Big(\sin(\sqrt{k}g(t))g''(t)+\sqrt{k}\cos(\sqrt{k}g(t))[(g'(t))^2-1]\Big)\geq0
\end{eqnarray}
on $[a,b]$. From this, we obtain, for all $t\in[a,b]$,
\begin{eqnarray}\label{orkpls}
  g''(t)\geq\sqrt{k}\Big(1-(g'(t))^2\Big)\cot(\sqrt{k}g(t))
\end{eqnarray}
that is (\ref{inprp22kpls}), as required.

We now consider the other direction of the case (i) of part 2 of the proposition.
If we have the differential inequality (\ref{inprp22kpls}) for all $t\in[a,b]$,
then by multiplication of both sides of (\ref{inprp22kpls}) with $\cos(\sqrt{k}t)\sin(\sqrt{k}g(t))>0$ (in the region we consider), we get
\begin{eqnarray}
  \cos(\sqrt{k}t)\Big(\sin(\sqrt{k}g(t))g''(t)+\sqrt{k}\big[(g'(t))^2-1\big]\cos(\sqrt{k}g(t))\Big)\geq0.
\end{eqnarray}
The left hand side of this inequality equals to the derivative of the function $H(t)$ defined in (\ref{Htkpls}).
Hence one has $H'(t)\geq0$ and consequently, the function $H(t)$ is non-decreasing. On the other hand,
the function $H(t)$ can be written as
\begin{eqnarray}
  H(t)=-\frac{\psi'(t)}{\phi'(t)},
\end{eqnarray}
where
\begin{eqnarray}\label{psiphikpls}
  \psi(t):=\frac{\cos(\sqrt{k}g(t))}{\cos(\sqrt{k}t)}\quad{\rm and}\quad\phi(t):=\tan(\sqrt{k}t).
\end{eqnarray}
Now, by using the extended mean value theorem on the intervals $[t_1,t]$ and $[t,t_2]$ for $t_1<t<t_2$, we get
\begin{eqnarray}\label{etakpls}
  \frac{\psi(t)-\psi(t_1)}{\phi(t)-\phi(t_1)}=\frac{\psi'(\eta)}{\phi'(\eta)}=-H(\eta),\quad{\rm for\,\,some}\quad\eta\in(t_1,t)
\end{eqnarray}
and
\begin{eqnarray}\label{xikpls}
  \frac{\psi(t_2)-\psi(t)}{\phi(t_2)-\phi(t)}=\frac{\psi'(\xi)}{\phi'(\xi)}=-H(\xi),\quad{\rm for\,\,some}\quad\xi\in(t,t_2).
\end{eqnarray}
Since the function $H(t)$ is non-decreasing, we have $-H(\eta)\geq-H(\xi)$, i.e.
\begin{eqnarray}
  \frac{\psi(t)-\psi(t_1)}{\phi(t)-\phi(t_1)}\geq\frac{\psi(t_2)-\psi(t)}{\phi(t_2)-\phi(t)}
\end{eqnarray}
for all $t\in(t_1,t_2)$, which we can also write in the form
\begin{eqnarray}\label{gvngthrsltkpls}
  \psi(t)\geq\frac{1}{\phi(t_2)-\phi(t_1)}\Big([\phi(t_2)-\phi(t)]\psi(t_1)+[\phi(t)-\phi(t_1)]\psi(t_2)\Big)
\end{eqnarray}
since the function $\phi(t)=\tan(\sqrt{k}t)$ is strictly increasing. Using the definitions of $\psi$
and $\phi$ in both sides of (\ref{gvngthrsltkpls}), we get
\begin{eqnarray}
  \frac{\cos(\sqrt{k}g(t))}{\cos(\sqrt{k}t)}\geq\frac{\mathcal{B}}{\tan(\sqrt{k}t_2)-\tan(\sqrt{k}t_1)},
\end{eqnarray}
where
\begin{eqnarray}
  \mathcal{B}&=&\Big(\tan(\sqrt{k}t_2)-\tan(\sqrt{k}t)\Big)\frac{\cos(\sqrt{k}g(t_1))}{\cos(\sqrt{k}t_1)}\nonumber\\
  &&+\Big(\tan(\sqrt{k}t)-\tan(\sqrt{k}t_1)\Big)\frac{\cos(\sqrt{k}g(t_2))}{\cos(\sqrt{k}t_2)}.
\end{eqnarray}
Here by using
\begin{eqnarray}
  \!\!\!\!\!\!\!\!\!\!\!\!\!\!\!\!&\,&g(t_1)=g^{t_1,t_2}_k(t_1)=\frac{1}{\sqrt{k}}\arccos\Big(\sqrt{k}\Big[u\cos(\sqrt{k}t_1)+v\sin(\sqrt{k}t_1)\Big]\Big)\\
  \!\!\!\!\!\!\!\!\!\!\!\!\!\!\!\!&&{\rm (and\,\,the\,\,same\,\,for}\,\,t_2{\rm )},\nonumber
\end{eqnarray}
we get
\begin{eqnarray}
  \frac{\cos(\sqrt{k}g(t))}{\cos(\sqrt{k}t)}\geq\sqrt{k}\Big(u+v\tan(\sqrt{k}t)\Big).
\end{eqnarray}
Since $\cos(\sqrt{k}t)>0$ (in the region we consider), we get
\begin{eqnarray}
  g(t)\leq\frac{1}{\sqrt{k}}\arccos\Big(\sqrt{k}\Big[u\cos(\sqrt{k}t)+v\sin(\sqrt{k}t)\Big]\Big)=g^{t_1,t_2}_k(t)
\end{eqnarray}
on $[t_1,t_2]$.

The proof of the case (ii) in the second part of the proposition
is proved along the same lines, only by reversing the inequalities.


\begin{thebibliography}{99}

\bibitem{BBI} Burago, D., Burago, Y., Ivanov, S.: A Course in Metric Geometry. American Mathematical Society (2001)

\bibitem{Papadopoulos} Papadopoulos, A.: Metric Spaces, Convexity and Nonpositive Curvature. European Mathematical Society (2014)


\end{thebibliography}
\end{document}